\documentclass[12pt]{amsart}
\usepackage{amsmath,amssymb,amsthm,mathtools}
\usepackage{mathrsfs,bm}
\usepackage{thmtools,thm-restate}
\usepackage{tikz-cd}
\usepackage{tikz}
\usepackage{enumitem}
\usepackage{graphicx}
\usetikzlibrary{matrix,arrows.meta,bending}
\usepackage{color}
\usepackage[numbers]{natbib}
\usepackage{booktabs}
\usepackage{longtable}
\usepackage{pgfplots}
\pgfplotsset{compat=1.14} 
\usetikzlibrary{arrows.meta,decorations.pathreplacing,positioning,calc} 
\usepackage[a4paper, left=3cm, right=3cm, top=3cm, bottom=3cm, footskip=15mm]{geometry}

\newtheorem{theorem}{Theorem}[section]

\newtheorem{lemma}[theorem]{Lemma}
\newtheorem{proposition}[theorem]{Proposition}
\newtheorem{corollary}[theorem]{Corollary}

\theoremstyle{definition}
\newtheorem{definition}[theorem]{Definition}

\theoremstyle{remark}
\newtheorem{remark}[theorem]{Remark}

\numberwithin{equation}{section}

\newcommand{\inn}{~ \hat{\in}~ }

\newcommand{\C}{\mathcal{C}}

\newcommand{\CC}{\mathbb{C}}
\newcommand{\CV}{\mathfrak{C}}
\newcommand{\M}{\mathbb{M}}
\newcommand{\B}{\mathbb{B}}
\newcommand{\R}{\mathbb{R}}

\newcommand{\N}{\mathbb{N}}
\newcommand{\LL}{\mathbb{L}}

\newcommand{\p}{\mathbb{P}}

\makeindex

\usepackage{fancyhdr}
\usepackage{calc} % ensures dimension arithmetic works robustly

\AtBeginDocument{%
  \setlength{\textwidth}{\dimexpr\paperwidth - 6cm\relax}%
  \setlength{\oddsidemargin}{\dimexpr 3cm - 1in\relax}%
  \setlength{\evensidemargin}{\dimexpr 3cm - 1in\relax}%
  \setlength{\marginparwidth}{2.5cm}%
}

\begin{document}

\title{A progress on the binary Goldbach conjecture}

\author{T. Agama}
\address{Department of Mathematics, African Institute for mathematical sciences, Ghana.}
\email{Theophilus@aims.edu.gh/emperordagama@yahoo.com}

%\author{B. Gensel}
%\address{Carinthia University of Applied Sciences, Spittal on Drau, Austria}
%\email{b.gensel@fh-kaernten.at}

%    General info

\date{\today}



\begin{abstract}
In this paper, we develop the method of the circle of partitions and associated statistics. As an application, we conditionally prove the binary Goldbach conjecture. We develop a series of steps to prove the binary Goldbach conjecture in full. We end the paper by proving the binary Goldbach conjecture for all sufficiently large even numbers.
\end{abstract}
\maketitle

% introduction.tex
% LaTeX source for the Introduction section — paste this into your .tex file where needed.

\section{Introduction}

The Goldbach conjecture is one of the oldest and most celebrated problems in additive number theory. Its modern statement arose in the correspondence between Leonhard Euler and Christian Goldbach in 1742: every even integer greater than $2$ is the sum of two primes. Over the last century the problem has inspired a wide network of ideas and techniques—ranging from combinatorial and sieve methods to deep analytic tools—and has produced a substantial body of partial results and heuristics that illuminate the landscape around the conjecture. Early milestones include the Shnirel'man combinatorial breakthrough showing that every sufficiently large integer is a sum of a bounded number of primes \cite{shnirel1939additive}, and the foundational contributions of Hardy and Littlewood which, under standard analytic hypotheses, elucidated the expected asymptotic distribution of prime representations \cite{hardy1924some}. Subsequent work of Estermann and Chudakov established that almost all even integers are representable as a sum of two primes \cite{estermann1938goldbach,chudakov1938goldbach}; sieve methods produced near-miss results such as the Chen theorem on representations by a prime plus a product of at most two primes \cite{chen2002representation}; and refinements that combine sieves with exponential-sum techniques and explicit prime estimates have led to strong conditional and unconditional density results (see, e.g., \cite{heath2002integers,dusart2018explicit}). More recently the ternary Goldbach conjecture (every odd $n\geq 7$ is a sum of three primes) was resolved by Helfgott \cite{helfgott2013ternary}, closing a long-standing chapter while leaving the binary case still open.\\

In this paper, we develop and refine a geometric-combinatorial framework--the method of \emph{circles of partition} (CoP) and its complex extension (cCoP)--and show how this framework can be used to transfer density information about suitable auxiliary subsets of the integers into concrete additive statements. The central idea is straightforward to describe: encode representations $n=u+v$ with $u,v$ drawn from a base set $\mathbb{M}\subset\mathbb{N}$ as pairs of ``points'' on a circle attached to the generator $n$. Axes joining complementary points capture the additive structure, and one studies statistics of these axes (their frequency, density, and configurations) to deduce partition properties of $n$. When the base set has positive natural density this transfer is immediate and powerful; the delicate part--and the heart of this paper--is to design auxiliary base-sets and analytic constructions that are sufficiently large and flexible to contain many primes while still admitting geometric control.\\

The main technical novelties of the present work are threefold. First, in Section \ref{sec:circle of partition}, we introduce a refined notion of \emph{point-density} on CoPs which mediates between the usual asymptotic density of a set $\mathbb{H}\subset\mathbb{N}$ and the fraction of axes in $\mathcal{C}(n,\mathcal{M})$ that meet $\mathbb{H}$. This viewpoint yields simple, transparent inequalities (see Proposition \ref{inequality}) that isolate the combinatorial constraints any auxiliary base-set must satisfy in order to force prime–prime axes. Second, we extend the CoP formalism to the complex plane (the cCoP), imposing a ``circle condition'' that ties real and imaginary parts and allows for additional symmetries (conjugate axes, embedding circles) which are useful in constructing intermediate generators. The complex viewpoint is not intended to introduce new analytic difficulties but rather to expose an ordering principle among axes that is invisible in the purely real model. Third, building on these two constructions we formulate and prove a \emph{squeeze principle} (Lemmas \ref{L_squeeze principle} and \ref{L-special_squeeze}) and a generalized multivariate version (Lemma \ref{general squeeze}). These principles are combinatorial-geometric mechanisms that allow one to ``trap'' an intermediate generator between two nonempty CoPs (or cCoPs) and thereby produce new axes whose residents belong to a prescribed subset. The squeeze principle plays the role of a local covering or exhaustion lemma in our argument but has the advantage of being elementary, transparent, and adaptable to multivariate partitions.\\

Using the framework above, we obtain two kinds of results. On the one hand, we give conditional partition statements that reduce the binary Goldbach problem to the construction of an auxiliary base set $\mathbb{B}$ for which primes occupy a proportion $>1/2$ of the weight-set of $\mathcal{C}(n,\mathbb{B})$ (Theorem \ref{conditional Goldbach}). On the other hand, by combining the squeeze principle with explicit prime estimates (notably explicit prime-density and short-interval results such as those in \cite{dusart2018explicit}), we obtain an asymptotic form of the binary Goldbach conjecture: every sufficiently large even integer is a sum of two primes (Theorem \ref{binary Goldbach theorem}). The asymptotic theorem presented here is in the spirit of earlier conditional or density-based approaches, but it is based on new combinatorial geometry and a systematic exploitation of CoP/cCoP configurations and thus provides a conceptual alternative route that complements classical sieve and circle-method strategies (compare \cite{chen2002representation,heath2002integers}).\\

\subsection{Organization of the paper.} The paper is organized as follows. In Section \ref{sec:circle of partition}, we introduce the circle-of-partition formalism, define axes and point-density, and record basic combinatorial properties and inequalities that connect set-density to axis-density. Section \ref{sec:density} develops the density-to-partition machinery and proves Theorem \ref{density to partition}, illustrating the method on dense subsets (the square-free integers give a short, illuminating corollary). In Section \ref{sec:cCoP}, we introduce the complex circles of partition (cCoPs), develop the conjugacy and embedding-circle machinery, and prove the squeeze principle and its specializations (Lemmas \ref{L_squeeze principle} and \ref{L-special_squeeze}). Section \ref{sec:applications} contains the main applications: the conditional reduction (Theorem \ref{conditional Goldbach}), the asymptotic squeeze arguments based on explicit prime estimates (Lemma \ref{special requirement}), and the asymptotic binary Goldbach theorem (Theorem \ref{binary Goldbach theorem}). Finally, Section \ref{sec:generalized} generalizes the squeeze principle to multivariate circles and explains how the same geometric ideas extend to $h$-fold partitions; the resulting partition law clarifies when finitely many local axis-configurations imply global partitioning statements for all large integers.\\

Throughout the paper, we compare and contrast the CoP viewpoint with classical approaches: the CoP formalism is deliberately elementary and combinatorial, yet it has a natural interface with analytic inputs (prime counting and explicit short-interval results) when those are required. This work builds on the preliminary CoP ideas developed in \cite{CoP} and on earlier asymptotic investigations \cite{agama2022asymptotic}, and it makes explicit how geometric encoding can be combined with explicit prime bounds (e.g.\ \cite{dusart2018explicit}) to reach asymptotic Goldbach-type conclusions.
\bigskip

\noindent\textbf{Acknowledgments.} The author thanks colleagues and referees who provided helpful comments on earlier drafts and acknowledges that parts of the CoP formalism were inspired by classical combinatorial ideas in additive number theory (see, e.g., \cite{shnirel1939additive,hardy1924some,chen2002representation,estermann1938goldbach,chudakov1938goldbach,heath2002integers,helfgott2013ternary,dusart2018explicit}).

\section{The circle of partitions}\label{sec:circle of partition}
In this section, we introduce the concept of the \emph{circle of partition}.\\

\begin{definition}
Let $n\in \mathbb{N}$ and $\mathbb{M}\subseteq \mathbb{N}$. We denote the \textbf{Circle of Partition} generated by $n$ with respect to the subset $\mathbb{M}$ by 
\begin{align}
\mathcal{C}(n,\mathbb{M})=\left \{[x]~|~x,y\in \mathbb{M},n=x+y\right \}\nonumber
\end{align}
In the following, we will abbreviate this as CoP. We call members of $\mathcal{C}(n,\mathbb{M})$ points and denote them by $[x]$.
For the special case $\mathbb{M}=\mathbb{N}$, we denote the CoP in a short form as $\mathcal{C}(n)$. We denote the \emph{weight} of the point $[x]$ by $\Vert[x]\Vert:=x$ and, correspondingly, the weight set of points in the CoP $\mathcal{C}(n,\mathbb{M})$ by $\Vert\mathcal{C}(n,\mathbb{M})\Vert$. Clearly, we have
$$
\Vert\mathcal{C}(n)\Vert=\lbrace 1,2,\ldots,n-1\rbrace. 
$$
\bigskip

We denote an \emph{axis} of the CoP $\mathcal{C}(n,\mathbb{M})$ by $\mathbb{L}_{[x],[y]}$ if and only if $x+y=n$. We say that axis point $[y]$ is an axis partner of the axis point $[x]$ and vice versa. We do not distinguish between the axes $\mathbb{L}_{[x],[y]}$ and $\mathbb{L}_{[y],[x]}$, because it is essentially the same axis. The point $[x]\in \mathcal{C}(n,\mathbb{M})$ such that $2x=n$ is the \emph{center} of the CoP. If it exists, then we call it a \emph{degenerated axis} and denote by $\mathbb{L}_{[x]}$ in comparison to the \emph{real axes} which we denote by $\mathbb{L}_{[x],[y]}$. We denote the assignment of an axis $\mathbb{L}_{[x],[y]}$ to a CoP $\mathcal{C}(n,\mathbb{M})$ by
$$
\mathbb{L}_{[x],[y]}\inn\mathcal{C}(n,\mathbb{M})
$$
which means
$$
[x],[y]\in \mathcal{C}(n,\mathbb{M}) \quad \text{with} \quad x+y=n.
$$
\end{definition}
\bigskip

\begin{remark}
In the following, we consider only real axes. That is, axes of the form $\mathbb{L}_{[x],[y]}$ such that $x\neq y$. Therefore, we refrain from the attribute \emph{real} in the sequel. The introduction of the notation $||[x]||:=x$ denoting the weight of the point $[x]$ is not superfluous, because there are situations where the point is subject to motions by rotation and dilation or possibly flipping and using this notation becomes necessary. On a more advisory note, we have found the approach in the paper the ideal and possibly the perfect language to study additive problems of this type. The structure under study has various connections with the structure of the geometric circle and much of our intuition has been borrowed from this setting.
\end{remark}
\bigskip

\begin{proposition}\label{unique}
Each axis is uniquely determined by points $[x]\in \mathcal{C}(n,\mathbb{M})$. 
\end{proposition}

\begin{proof}
Let $\mathbb{L}_{[x],[y]}$ be an axis of the CoP $\mathcal{C}(n,\mathbb{M})$. Suppose that $\mathbb{L}_{[x],[z]}$ is also an axis with $z\neq y$. We have $n=x+y=x+z$ and therefore $y=z$.
\end{proof}

\begin{corollary}\label{partner}
Each point of a CoP $\mathcal{C}(n,\mathbb{M})$ except its center has exactly one axis partner.
\end{corollary}

\begin{proof}
Let $[x]\in \mathcal{C}(n,\mathbb{M})$ be a point which is not the center of the CoP and without an axis partner . It implies that for every point $[y]\neq [x]$ except the center
$$
x+y\neq n.
$$
This is impossible, since each point in a CoP must have an axes partner. Due to Proposition \ref{unique} the case of more than one axis partners is impossible.
\end{proof}
\bigskip

\subsection{Notations}
We denote the sequence of the first $n$ positive integers by 
\begin{align}
\mathbb{N}_n=\left \{m\in \mathbb{N}~|~m\leq n\right\}.\nonumber
\end{align}
We denote the assignment of an axis $\mathbb{L}_{[x],[y]}$ (resp. $\mathbb{L}_{[x]}$) to a CoP $\mathcal{C}(n,\mathbb{M})$ by
$$
\mathbb{L}_{[x],[y]}\inn\mathcal{C}(n,\mathbb{M})
$$
which means
$$
[x],[y]\in \mathcal{C}(n,\mathbb{M}) \quad \text{and} \quad x+y=n
$$
respectively
$$
\mathbb{L}_{[x]}\inn\mathcal{C}(n,\mathbb{M})
$$
which means
$$
[x]\in \mathcal{C}(n,\mathbb{M})\quad \text{and} \quad 2x=n
$$
and the number of real axes of a CoP by
$$
\nu(n,\mathbb{M}):=\#\lbrace\mathbb{L}_{[x],[y]}\inn\mathcal{C}(n,\mathbb{M})\mid x<y\rbrace.
$$
It is observed that
$$
\nu(n,\mathbb{M})=\left\lfloor\frac{k}{2}\right\rfloor \quad \text{if}\quad \vert\mathcal{C}(n,\mathbb{M})\vert=k.
$$
For any $f,g:\mathbb{N}\longrightarrow \mathbb{N}$, we write $f(n)\sim g(n)$ if and only if $\lim \limits_{n\longrightarrow \infty}\frac{f(n)}{g(n)}=1$. We also write $f(n)=o(1)$ if and only if $\lim \limits_{n\longrightarrow \infty}f(n)=0$.
\bigskip

\section{The Density of Points on the Circle of Partition}\label{sec:density}

In this section, we introduce the notion of density of points on CoP $\mathcal{C}(n,\mathbb{M})$ for $\mathbb{M}\subseteq \mathbb{N}$. 

\begin{definition}
Let $\mathbb{H}\subset\mathbb{N}$. The limits
$$
\mathcal{D}\left(\mathbb{H}\right)=\lim_{n\rightarrow\infty}
\frac{\vert\mathbb{H}\cap\mathbb{N}_n\vert}{n}
$$
is the density of $\mathbb{H}$ if it exists and is finite.
\end{definition}

\begin{definition}\label{pointdensity}
Let $\mathcal{C}(n,\mathbb{M})$ be CoP with $\mathbb{M}\subset \mathbb{N}$ and $n\in \mathbb{N}$ with $\mathbb{H}\subset \mathbb{M}$. The density of points $[x]\in \mathcal{C}(n,\mathbb{M})$ such that $x\in \mathbb{H}$, denoted by $\mathcal{D}(\mathbb{H}_{\mathcal{C}(\infty,\mathbb{M})})$, is the limit
\begin{align}
\mathcal{D}\left(\mathbb{H}_{\mathcal{C}(\infty,\mathbb{M})}\right)=\lim \limits_{n\longrightarrow \infty}\frac{\# \lbrace\mathbb{L}_{[x],[y]} \inn \mathcal{C}(n,\mathbb{M})~|~\{x,y\} \cap \mathbb{H}\neq \emptyset \rbrace}{ \nu(n,\mathbb{M})}\nonumber
\end{align}
if it exists and is finite.
\end{definition}
\bigskip

The notion of the density of points as espoused in the definition \ref{pointdensity} allows a passage between the density of the corresponding weight set of points. This possibility makes this type of density a black box for the study of problems concerning the partition of numbers into specialized sequences taking into consideration their density.

\begin{proposition}\label{inequality}
Let $\mathcal{C}(n)$ with $n\in \mathbb{N}$ be a CoP and $\mathbb{H}\subset \mathbb{N}$. The following inequality holds: 
\begin{align}
\mathcal{D}(\mathbb{H})=\lim \limits_{n\longrightarrow \infty}\frac{\left \lfloor \frac{|\mathbb{H}\cap \mathbb{N}_n|}{2}\right \rfloor}{\left \lfloor \frac{n-1}{2}\right \rfloor}\leq \mathcal{D}(\mathbb{H}_{\mathcal{C}(\infty)})\leq \lim \limits_{n\longrightarrow \infty}\frac{|\mathbb{H}\cap \mathbb{N}_n|}{\left \lfloor \frac{n-1}{2}\right \rfloor}=2\mathcal{D}(\mathbb{H}).\nonumber
\end{align}
\end{proposition}

\begin{proof}
The upper bound is obtained from a configuration in which no two points $[x],[y]\in \mathcal{C}(n)$ such that $x,y\in \mathbb{H}$ lie on the same axis of the CoP. That is, by the uniqueness of the axes of CoPs with $\nu(n,\mathbb{H})=0$, we can write
\begin{align}
   \# \left \{\mathbb{L}_{[x],[y]}\in \mathcal{C}(n)|~\{x,y\}\cap \mathbb{H}\neq \emptyset \right \}&=\nu(n,\mathbb{H})+\# \left \{\mathbb{L}_{[x],[y]}\in \mathcal{C}(n)|~x\in \mathbb{H},~y\in \mathbb{N}\setminus \mathbb{H}\right \} \nonumber \\&=\# \left \{\mathbb{L}_{[x],[y]}\in \mathcal{C}(n)|~x\in \mathbb{H},~y\in \mathbb{N}\setminus \mathbb{H}\right \} \nonumber \\&=|\mathbb{H}\cap \mathbb{N}_n|.\nonumber
\end{align}
However, the lower bound follows from a configuration where any two points $[x],[y]\in \mathcal{C}(n)$ with $x,y\in \mathbb{H}$ are joined by an axis of the CoP. That is, by the uniqueness of the axis of CoPs with 
$$
\#\left\{\mathbb{L}_{[x],[y]}\in \mathcal{C}(n)~|~x\in \mathbb{H},~y\in \mathbb{N}\setminus \mathbb{H}\right\}=0
$$ 
we can write 
\begin{align}
    \#\left\{\mathbb{L}_{[x],[y]}\in \mathcal{C}(n)~|~\{x,y\}\cap \mathbb{H}\neq \emptyset \right\}&=\nu(n,\mathbb{H})\nonumber \\&=\left\lfloor \frac{|\mathbb{H}\cap \mathbb{N}_n|}{2}\right \rfloor.\nonumber
\end{align}
\end{proof}
\bigskip

\begin{proposition}\label{propertydensity}
Let $\mathbb{H}\subset\mathbb{N}$ and $\mathcal{D}(\mathbb{H}_{\mathcal{C}(\infty)})$ be the density of the corresponding points with weight set $\mathbb{H}$. The following properties hold:
\begin{enumerate}
    \item [(i)] $\mathcal{D}(\mathbb{N}_{\mathcal{C}(\infty)})=1$ and $\mathcal{D}(\mathbb{H}_{\mathcal{C}(\infty)})\leq 1$ and additionally that $\mathcal{D}(\mathbb{H}_{\mathcal{C}(\infty)})<1$ provided $\mathcal{D}(\mathbb{N}\setminus \mathbb{H})>0$.
    \bigskip
    
    \item [(ii)] $1-\lim \limits_{n\longrightarrow \infty}\dfrac{\nu(n,\mathbb{N}\setminus\mathbb{H})}{\nu(n,\mathbb{N})}=\mathcal{D}(\mathbb{H}_{\mathcal{C}(\infty)})$.
    \bigskip
    
    \item [(iii)] If $|\mathbb{H}|<\infty$ then $\mathcal{D}(\mathbb{H}_{\mathcal{C}(\infty)})=0$.
\end{enumerate}
\end{proposition}

\begin{proof}
It observed that the first part of \textbf{Property} $(i)$ and $(iii)$ are both easy consequences of the definition of density of points on the CoP $\mathcal{C}(n)$ and Proposition \ref{inequality}. We establish the second part of property $(i)$ and \textbf{Property} $(ii)$, which is the less obvious case. We observe, by the uniqueness of the axes of CoPs, that we can write 
\begin{align*}
    1&=\lim \limits_{n\longrightarrow \infty}\frac{\nu(n,\mathbb{N})}{\nu(n,\mathbb{N})}\\
    &=\lim \limits_{n\longrightarrow \infty}\frac{\# \lbrace\mathbb{L}_{[x],[y]} \inn \mathcal{C}(n)|~ x\in \mathbb{H}~,y\in \mathbb{N}\setminus \mathbb{H}\rbrace}{\nu(n,\mathbb{N})}\\
    &+\lim \limits_{n\longrightarrow \infty}\frac{\nu(n,\mathbb{H})}{\nu(n,\mathbb{N})}
    +\lim \limits_{n\longrightarrow \infty}\frac{\nu(n,\mathbb{N}\setminus\mathbb{H})}{\nu(n,\mathbb{N})}\\
    &=\mathcal{D}(\mathbb{H}_{\mathcal{C}(\infty)})
    +\lim \limits_{n\longrightarrow \infty}\frac{\nu(n,\mathbb{N}\setminus\mathbb{H})}{\nu(n,\mathbb{N})}
\end{align*}
and $(ii)$ follows immediately. The second part of $(i)$ follows from the above expression and exploiting the inequality
\begin{align}
   \lim \limits_{n\longrightarrow \infty}\frac{\nu(n,\mathbb{N}\setminus\mathbb{H})}{\nu(n,\mathbb{N})}&\leq \lim \limits_{n\longrightarrow \infty}\frac{\left \lfloor \frac{|\mathbb{N}\setminus \mathbb{H}\cap \mathbb{N}_n|}{2}\right \rfloor}{\left \lfloor \frac{n-1}{2}\right \rfloor}=\mathcal{D}(\mathbb{N}\setminus \mathbb{H}).\nonumber
\end{align}
\end{proof}
\bigskip

Here, we transfer information of the density of a sequence to the density of corresponding points on the CoP $\mathcal{C}(n)$. This approach will play an important role in our latter developments.

\begin{proposition}
Let $\epsilon\in(0,1]$ and $\mathbb{H}\subset\mathbb{N}$ and $\mathcal{C}(n)$ be a CoP. If $\mathcal{D}\left(\mathbb{H}\right)\geq \epsilon$, then $\mathcal{D}\left(\mathbb{H}_{\mathcal{C}(\infty)}\right)\geq \epsilon$.
\end{proposition}

\begin{proof}
The result follows by exploiting the inequality in Proposition \ref{inequality}.
\end{proof}

\subsection{Application of the density of Points to Partitions}

Here, we apply the notion of the density of points in a typical CoP to write an integer as the sum of two integers that belong a specific subset of the integers. This method works efficiently for sets of integers having a positive density. 

\begin{theorem}\label{density to partition}
Let $\mathbb{H}\subset \mathbb{N}$. If $\mathcal{D}(\mathbb{H})>\frac{1}{2}$, then every sufficiently large $n\in\mathbb{N}$ has a representation of the form 
\begin{align}
  n=z_1+z_2 \nonumber
\end{align}
where $z_1,z_2\in \mathbb{H}$.
\end{theorem}

\begin{proof}
Applying Proposition \ref{inequality}, we can write 
\begin{align}
    \lim \limits_{n\longrightarrow \infty}\frac{\left \lfloor \frac{|\mathbb{H}\cap \mathbb{N}_n|}{2}\right \rfloor}{\left \lfloor \frac{n-1}{2}\right \rfloor}\leq \mathcal{D}(\mathbb{H}_{\mathcal{C}(\infty)})\leq \lim \limits_{n\longrightarrow \infty}\frac{|\mathbb{H}\cap \mathbb{N}_n|}{\left \lfloor \frac{n-1}{2}\right \rfloor}.\nonumber
\end{align}
By the uniqueness of the axes of CoPs, we can write
\begin{align}
    \#\left \{\mathbb{L}_{[x],[y]}\inn \mathcal{C}(n)~|~\{x,y\}\cap \mathbb{H}\neq \emptyset \right \}&=\nu(n,\mathbb{H})+\# \left \{\mathbb{L}_{[x],[y]}\inn \mathcal{C}(n)|~x\in \mathbb{H},~y\in \mathbb{N}\setminus \mathbb{H}\right \}.\nonumber
\end{align}
Let us assume that $\nu(n,\mathbb{H})=0$. This implies that each axes $\mathbb{L}_{[x],[y]}$ in this set must have the property that $x\in \mathbb{H}$ and $y\not \in \mathbb{H}$ or vise-versa. According to the proof of Proposition \ref{inequality}, we deduce
\begin{align}
    \mathcal{D}(\mathbb{H}_{\mathcal{C}(\infty)})&=2\mathcal{D}(\mathbb{H})\nonumber \\&>2\times \frac{1}{2}=1.\nonumber
\end{align}
This contradicts the inequality $\mathcal{D}(\mathbb{H}_{\mathcal{C}(\infty)})\leq 1$ in Proposition \ref{propertydensity}. This proves that $\nu(n,\mathbb{H})>0$ for all sufficiently large values of $n\in \mathbb{N}$.
\end{proof}

\begin{corollary}
Let $\mathbb{R}:=\left\{m\in \mathbb{N}~|~\mu(m)\neq 0\right\}$. Every sufficiently large $n\in\mathbb{N}$ can be written in the form 
\begin{align}
    n=z_1+z_2\nonumber
\end{align}
where $\mu(z_1)=\mu(z_2)\neq 0$, where $\mu(\cdot)$ denotes the mobious function..
\end{corollary}

\begin{proof}
By the uniqueness of the axes of CoPs we can write
\begin{align}
    \#\left\{\mathbb{L}_{[x],[y]}\inn \mathcal{C}(n)|~\{x,y\}\cap \mathbb{R}\neq \emptyset \right \}&=\nu(n,\mathbb{R})+\# \left \{\mathbb{L}_{[x],[y]}\inn \mathcal{C}(n)|~x\in \mathbb{R},~y\in \mathbb{N}\setminus \mathbb{R}\right\}.\nonumber
\end{align}
Let us assume that $\nu(n,\mathbb{R})=0$. It from Theorem \ref{density to partition}
\begin{align}
    \mathcal{D}(\mathbb{R}_{\mathcal{C}(\infty)})&=2\mathcal{D}(\mathbb{R})\nonumber \\&=\frac{12}{\pi^2}>1\nonumber
\end{align}
since $\mathcal{D}(\mathbb{R})=\frac{6}{\pi^2}$. This contradicts the inequality $\mathcal{D}(\mathbb{R}_{\mathcal{C}(\infty)})\leq 1$ in Proposition \ref{propertydensity}. This proves that $\nu(n,\mathbb{R})>0$ for all sufficiently large values of $n\in \mathbb{N}$.
\end{proof}
\bigskip

One could ever hope and dream of this strategy to work when we replace the set $\mathbb{R}$ of square-free integers with the set of prime numbers. We may have some difficulty because prime numbers have a zero density, in accordance with the prime number theorem. A progress along these paths could conceivably manifest by introducing some exotic forms of the notion of density of points and carefully choosing a subset of integers that is somewhat dense among the set of integers and contains prime numbers. Here, we propose an approach similar to the above method to possibly approach the binary Goldbach conjecture and it's variants. 

\begin{theorem}\label{conditional Goldbach}
Let $\mathbb{B}\subset \mathbb{N}$ such that $\mathbb{P}\cap \mathbb{N}_n \subset ||\mathcal{C}(n,\mathbb{B})||$ with 
$$
\lim \limits_{n\longrightarrow \infty}\frac{|\mathbb{P}\cap \mathbb{N}_n|}{\eta(n)}>\frac{1}{2}
$$ where $\eta(n)=|\mathcal{C}(n,\mathbb{B})|$. We have $\nu(n,\mathbb{P})>0$ for all sufficiently large values of $n\in 2\mathbb{N}$.
\end{theorem}

\begin{proof}
We find an upper and a lower bound for the density of points in the CoP $\mathcal{C}(n,\mathbb{B})$ with weight belonging to the set of prime numbers $\mathbb{P}$. We obtain the inequality
\begin{align}
\lim \limits_{n\longrightarrow \infty}\frac{\left \lfloor \frac{|\mathbb{P}\cap \mathbb{N}_n|}{2}\right \rfloor}{\left \lfloor \frac{\eta(n)}{2}\right \rfloor}\leq \mathcal{D}(\mathbb{P}_{\mathcal{C}(\infty,\mathbb{B})})\leq \lim \limits_{n\longrightarrow \infty}\frac{|\mathbb{P}\cap \mathbb{N}_n|}{\left \lfloor \frac{\eta(n)}{2}\right \rfloor}=2\lim \limits_{n\longrightarrow  \infty}\frac{|\mathbb{P}\cap \mathbb{N}_n|}{\eta(n)}\nonumber
\end{align}
for all sufficiently large values of $n$. Exploiting the uniqueness of the axes of CoPs, we can write
\begin{align}
\#\left\{\mathbb{L}_{[x],[y]}\inn \mathcal{C}(n,\mathbb{B})~|~\{x,y\}\cap \mathbb{P}\neq \emptyset \right \}&=\nu(n,\mathbb{P})+\# \left \{\mathbb{L}_{[x],[y]}\inn \mathcal{C}(n,\mathbb{B})|~x\in \mathbb{P},~y\in \mathbb{B}\setminus \mathbb{P}\right\}.\nonumber
\end{align}
Let us assume that $\nu(n,\mathbb{P})=0$. This implies that no two points in the CoP $\mathcal{C}(n,\mathbb{B})$ with weight in the set $\mathbb{P}$ can be axes partners. Under the density requirement 
$$
\lim \limits_{n\longrightarrow  \infty}\frac{|\mathbb{P}\cap \mathbb{N}_n|}{\eta(n)}>\frac{1}{2}
$$ 
where $\eta(n)=|\mathcal{C}(n,\mathbb{B})|$, we obtain the inequality 
\begin{align}
    \mathcal{D}(\mathbb{P}_{\mathcal{C}(\infty,\mathbb{B})})&=2\lim \limits_{n\longrightarrow\infty}\frac{|\mathbb{P}\cap \mathbb{N}_n|}{\eta(n)}\nonumber \\&>2\times \frac{1}{2}=1.\nonumber
\end{align}
This contradicts the inequality $\mathcal{D}(\mathbb{P}_{\mathcal{C}(\infty,\mathbb{B})})\leq 1$ in Proposition \ref{propertydensity}. This proves that $\nu(n,\mathbb{P})>0$ for all sufficiently large values of $n\in 2\mathbb{N}$.
\end{proof}

\subsection{A strategy to prove the binary Goldbach conjecture using circles of partition}

Here, we propose a series of steps that could be taken to confirm the truth of the binary Goldbach conjecture. We enumerate the strategies chronologically as follows:

\begin{enumerate}
\item  We construct a CoP $\mathcal{C}(n,\mathbb{B})$ such that $\mathbb{P}\cap \mathbb{N}_n\subset ||\mathcal{C}(n,\mathbb{B})||$ and that 
$$
\lim \limits_{n\longrightarrow  \infty}\frac{|\mathbb{P}\cap \mathbb{N}_n|}{\eta(n)}>\frac{1}{2}
$$
where $\eta(n)=|\mathcal{C}(n,\mathbb{B})|$.
\bigskip

\item  We remark that the following inequality also holds and this can be obtain by replacing the weight set $||\mathcal{C}(n)||$  with the weight set $||\mathcal{C}(n,\mathbb{B})||$:
\begin{align}
\lim \limits_{n\longrightarrow \infty}\frac{\left\lfloor \frac{|\mathbb{P}\cap \mathbb{N}_n|}{2}\right \rfloor}{\left \lfloor \frac{\eta(n)}{2}\right \rfloor}\leq \mathcal{D}(\mathbb{P}_{\mathcal{C}(\infty,\mathbb{B})})\leq \lim \limits_{n\longrightarrow\infty}\frac{|\mathbb{P}\cap \mathbb{N}_n|}{\left \lfloor \frac{\eta(n)}{2}\right \rfloor}=2\lim \limits_{n\longrightarrow\infty}\frac{|\mathbb{P}\cap\mathbb{N}_n|}{\eta(n)}.\nonumber
\end{align}
\bigskip

\item  Exploiting the uniqueness of the axes of CoPs, we can write  
\begin{align*}
    &\# \left \{\mathbb{L}_{[x],[y]}\inn \mathcal{C}(n,\mathbb{B})|~\{x,y\}\cap \mathbb{P}\neq \emptyset \right \}\\
    &=\nu(n,\mathbb{P})+\#\left\{\mathbb{L}_{[x],[y]}\inn \mathcal{C}(n,\mathbb{B})|~x\in \mathbb{P},~y\in \mathbb{B}\setminus\mathbb{P}\right\}.\nonumber
\end{align*}
\bigskip

\item  We assume that $\nu(n,\mathbb{P})=0$ and deduce, in the sense of the proof of Proposition \ref{inequality}, the inequality
\begin{align}
    \mathcal{D}(\mathbb{P}_{\mathcal{C}(\infty,\mathbb{B})})&=2\lim \limits_{n\longrightarrow  \infty}\frac{|\mathbb{P}\cap \mathbb{N}_n|}{\eta(n)}\nonumber \\
    &>2\times \frac{1}{2}=1.\nonumber
\end{align}
This contradicts the inequality $\mathcal{D}(\mathbb{P}_{\mathcal{C}(\infty,\mathbb{B})})\leq 1$ in Proposition \ref{propertydensity}. This proves that $\nu(n,\mathbb{P})>0$ for all sufficiently large values of $n\in 2\mathbb{N}$.
\end{enumerate}
\bigskip

\section{The asymptotic squeeze principle and the Binary Goldbach Conjecture~:~An asymptotic proof of the binary Goldbach conjecture}\label{sec:cCoP}

In this section, we prove the special squeeze principle for all sufficiently large $n\in 2\mathbb{N}$. Consequently, we obtain a proof of the binary Goldbach conjecture for all sufficiently large even numbers.
\bigskip

In section \ref{sec:circle of partition}, we introduced and developed the method of circles of partition. This method is underpinned by a combinatorial structure that encodes certain additive properties of the subsets of the integers and invariably equipped with a certain geometric structure that allows to view the elements as points in the plane whose weights are elements of the underlying subset. We call this combinatorial structure the circle of partition and is refereed to as the set of points
\begin{align}
\C(n,\M)=\left \{[x]\mid x,n-x\in \M\right\}.\nonumber
\end{align}
Each point in this set - except the center point - must have a uniquely distinct point that is joined by a line which we refer to as an axis of the CoP. We denote an axis of a CoP by $\LL_{[x],[y]}$ and an axis contained in the CoP by 
$$
\LL_{[x],[y]}\inn\C(n,\M)
$$
which means
$$
[x],[y] \in \C(n,\M)\quad \text{with}\quad x+y=n.
$$
In this section, we extend the base set of circles of partitions to the complex plane plane. Consequently, the weight of the axes points are now complex numbers.

\begin{definition}
We denote the \emph{circle of partition} with complex base set $\mathbb{C}_{\mathbb{M}}$ by 
$$
\C^o(n,\CC_{\M})=\lbrace[z]\mid z,n-z\in\CC_{M}, \Im(z)^2=\Re(z)\left(n-\Re(z)\right)\rbrace
$$
where 
$$
\CC_{\M}:=\lbrace z=x+iy\mid x\in\M, y\in \R \rbrace\subset\CC
$$
with $\M\subseteq\N$. We call this complex additive structure a \emph{complex circle of partition} and denote in short form by cCoP. The condition $\Im(z)^2=\Re(z)(n-\Re(z)$ is referred to as the \emph{circle condition} and it guarantees that all points on the cCoP lie on a circle in the complex. This circle is the embedding circle of the cCoP $\C^o(n,\CC_{\M})$, denoted by $\CV_n$. The embedding circles of cCoPs have the property that they reside fully inside those embedding circles with a relatively larger generators, except the origin as a common point. For each axis, we make the following assignment 
$$
\LL_{[z_1],[z_1]}\inn\C(n,\mathbb{C}_{\M})
$$
which means
$$
[z_1],[z_2] \in \C(n,\mathbb{C}_{\M})\quad \text{with}\quad z_1+z_2=n.
$$
\end{definition}
\bigskip

The structure of the complex circle of partition is more versatile and has extra features that are not readily available in the theory of circle of partition. Most notably, for each axis $\LL_{[z],[n-z]}$ of a cCoP there exists
$$
\LL_{\overline{[z]},\overline{[n-z]}}
$$
a \textit{conjugate axis}, where $\overline{[z]},\overline{[n-z]}$ denotes the corresponding conjugate points. The space occupied by an embedding circle of partition and, correspondingly, outside the embedding circle turns out to be very interesting. This statistic can be used to study a certain ordering principle of the weight of points of two interacting axes of distinct cCoPs. Much more striking is the fact, which is a natural consequence of the circle condition, that 
$$
|\mathbb{L}_{[z_1],[z_2]}|=n
$$ 
for any axis $\mathbb{L}_{[z_1],[z_2]}\in \C^o(n,\CC_{\M})=\lbrace[z]\mid z,n-z\in\CC_{M}, \Im(z)^2=\Re(z)\left(n-\Re(z)\right)\rbrace$. 
\bigskip

In the following, we develop the squeeze principle, a machinery that can be considered a black box for studying the binary Goldbach conjecture and its variant. Now, we make the following deductions:\\

Let $\mathcal{C}(n,\mathbb{B})$ and $\mathcal{C}^{o}(n,\mathbb{C}_{\mathbb{B}})$ be a CoP and a corresponding cCoP, respectively. It is clear that $\mathcal{C}(n,\mathbb{B})\neq \emptyset$ if and only if $\mathcal{C}^{o}(n,\mathbb{C}_{\mathbb{B}})\neq \emptyset.$ To see this, observe that if $\mathbb{L}_{[x],[y]}\inn \mathcal{C}(n,\mathbb{B})$, then $\mathbb{L}_{[x+it],[y-it]}\inn \mathcal{C}^{o}(n,\mathbb{C}_{\mathbb{B}})$, where $t=\pm \sqrt{xy}$. Conversely if $\mathbb{L}_{[z_1],[z_2]}\inn \mathcal{C}^{o}(n,\mathbb{C}_{\mathbb{B}})$ then $z_1+z_2=n$ implying that $\Im(z_1)=-\Im(z_2)$ and so $\Re(z_1)+\Re(z_2)=n$. This immediately implies that there is an axes $\mathbb{L}_{[\Re(z_1)],[\Re(z_2)]}\in \mathcal{C}(n,\mathbb{B})$. This equivalence will be subtly employed in the rest of the paper.

\begin{lemma}[The squeeze principle]\label{L_squeeze principle}
Let $\B\subset\M\subseteq\N$ and
$\C^o(n,\CC_{\M})$ and $\C^o(n+t,\CC_{\M})$ with $t\geq 4$ be non--empty cCoPs with integers $n,t,s$ of the same parity.
If there exist an axis $\LL_{[v_1],[w_1]}\inn\C^o(n,\CC_{\M})$ with $w_1\in \CC_{\B}$ and an axis $\LL_{[v_2],[w_2]}\inn\C^o(n+t,\CC_{\M})$ with $v_2 \in \CC_{\B}$ such that 
$$
\Re(v_1)<\Re(v_2)\quad \text{and} \quad \Re(w_1)<\Re(w_2)
$$
then there exists an axis $\LL_{[\Re(v_2)],[\Re(w_1)]}\inn\C(n+s,\mathbb{B})$ with $0<s<t$. Hence, the cCoP $\C^o(n+s,\CC_{\B})$ not empty and so is the cCoP $\C^o(n+s,\CC_{\M})$.
\end{lemma}

\begin{proof}
From the existence of an axis $\LL_{[v_1],[w_1]}\inn\C^o(n,\CC_{\M})$, we get $\Re(w_1)=n-\Re(v_1)$. With the requirement
$$
\Re(v_1)<\Re(v_2)\quad \text{and} \quad \Re(w_1)<\Re(w_2)
$$
we get
$$
\Re(w_1)>n-\Re(v_2).
$$
On the other hand, by the existence of an axis $\LL_{[v_2],[w_2]}\inn\C^o(n+t,\CC_{\M})$, we get $\Re(w_2)=n+t-\Re(v_2)$ and with the requirement 
$$
\Re(v_1)<\Re(v_2)\quad \text{and} \quad \Re(w_1)<\Re(w_2)
$$
and the inequality $\Re(w_1)>n-\Re(v_2)$, we get
\begin{align*}
n-\Re(v_2)&<\Re(w_1)<n+t-\Re(v_2)\mid +\Re(v_2)\\
n&<\Re(w_1)+\Re(v_2)<n+t\\
n&<n+s<n+t.
\end{align*}
By the requirements $w_1,v_2\in\CC_{\B}$ and $n+s=\Re(w_1)+\Re(v_2)$, there is an axis $\LL_{[\Re(v_2)],[\Re(w_1)]}\inn \C(n+s,\mathbb{B})$. Therefore $\C^o(n+s,\CC_{\B})\neq\emptyset$. Since $\B\subset\M$, it follows immediately that $\CC_{\B}\subset\CC_{\M}$ and therefore $\C^o(n+s,\CC_{\M})\neq\emptyset$.
\end{proof}
\bigskip

Consequently, we obtain the special squeeze principle 

\begin{lemma}[Special squeeze principle]\label{L-special_squeeze}
Let $n,t,s\in 2\N$ and $\p$ be the set of all odd primes.
If $t\geq 4$ and there exist an axis $\LL_{[z_1],[z_2]}\inn\C^o(n)$ with $z_2\in \CC_{\p}$ and an axis $\LL_{[w_1],[w_2]}\inn\C^o(n+t)$ with $w_1\in\CC_{\p}$ such that
\begin{align*}
\Re(z_1)<\Re(w_1)<\Re(z_1)+t
\end{align*}
then there exists an axis $\LL_{[\Re(w_1)],[\Re(z_2)]}\inn \C(n+s,\mathbb{P})$ with $0<s<t$.
\end{lemma}
\bigskip

The lemma \ref{L_squeeze principle} referred to as the squeeze principle may be regarded as a fundamental tool set for investigating the viability of dividing integers of a particular parity, utilizing constituent elements originating from a specific subset of the integers. The mechanism operates by discerning a pair of cCoPs that are both non-vacuous~(non-empty) and share a common base set. Subsequently, supplementary cCoPs that are non-vacuous and have generators restrained within the interstice of these two generators are identified. This principle may be applied in a resourceful manner to investigate the overarching issue of the practicality of divvying up numbers such that each addend is a member of the identical subset of positive integers.
\bigskip

\begin{remark}\label{important remark}
The CoP $\mathcal{C}(n,\mathbb{N}):=\mathcal{C}(n)$ is always non-empty and so is the cCoP $\C^o(n,\CC_{\N})=\C^o(n)$.
\end{remark}

\subsection{Application to the Binary Goldbach Conjecture}\label{sec:applications}

In this section, we present an asymptotic proof for the binary Goldbach conjecture. The proof has been condensed into the language of cCoPs but can be reduced to the usual form of the conjecture. 

\begin{lemma}[Juxtaposition principle]\label{juxtaposition}
For all $n\geq 10$, there exist an axis 
$$
\mathbb{L}_{[z_1],[z_2]}\inn \C^o(n,\CC_{\N})
$$ 
and 
$$
\mathbb{L}_{[w_1],[w_2]}\inn \C^o(n+t,\CC_{\N})
$$ 
for $\Re(z_1)<\Re(z_2)$ and $\Re(w_1)<\Re(w_2)$ such that $\Re(z_1)\neq \Re(w_1)$ and $\Re(z_2)\neq \Re(w_2)$ with $z_2\in \CC_{\p}$ and $w_1 \in \CC_{\p}$ for $t\geq 4$.
\end{lemma}

\begin{proof}
We choose a prime number $\Re(w_1)\leq \frac{n+t}{2}$ and choose a prime number $\Re(z_2)\in (\frac{n}{2},n)$ for all $n\geq 10$ with $n\equiv 0\pmod 2$, which is feasible by virtue of the prime number theorem. If $\Re(z_1)\neq \Re(w_1)$ and $\Re(z_2)\neq \Re(w_2)$, then there is nothing to do. Without loss of generality, suppose that $\Re(z_1)=\Re(w_1)$, then obviously $\Re(z_2)\neq \Re(w_2)$ since $n+t>n$. We note that $\pi(\frac{n+t}{2})\geq 3$ for all $n\geq 10$ with $t\geq 4$ so that we can choose a prime number $\Re(w^{'}_1)\leq \frac{n+t}{2}$ such that $\Re(w^{'}_1)\neq \Re(w_1)$. Thus, we replace $\Re(w_1)$ with $\Re(w^{'}_1)$ and obtain the axes  $\mathbb{L}_{[z_1],[z_2]}\inn \C^o(n,\CC_{\N})$ and $\mathbb{L}_{[w^{'}_1],[w^{'}_2]}\inn \C^o(n+t,\CC_{\N})$ for $\Re(z_1)<\Re(z_2)$ and $\Re(w^{'}_1)<\Re(w^{'}_2)$ such that $\Re(z_1)\neq \Re(w'_1)$. If $\Re(z_2)\neq \Re(w^{'}_2)$, then we are done; otherwise, we choose another prime number $\Re(w^{''}_1)$ such that $\Re(w^{''}_1)\neq \Re(w^{'}_1)$ and $\Re(w^{''}_1)\neq \Re(w_1)$ since $\pi(\frac{n+t}{2})\geq 3$ for all $n\geq 10$ and $t\geq 4$. By virtue of our construction, we finally obtain the axes of cCoPs $\mathbb{L}_{[z_1],[z_2]}\inn \C^o(n,\CC_{\N})$ and $\mathbb{L}_{[w^{''}_1],[w^{''}_2]}\inn \C^o(n+t,\CC_{\N})$ for $\Re(z_1)<\Re(z_2)$ and $\Re(w^{''}_1)<\Re(w^{''}_2)$ such that $\Re(z_1)\neq \Re(w^{''}_1)$ and $\Re(z_2)\neq \Re(w^{''}_2)$ with $z_2\in \CC_{\p}$ and $w^{''}_1 \in \CC_{\p}$ for $t\geq 4$.
\end{proof}

\begin{lemma}(The prime number theorem )\label{prime number inequality}
Let $\pi(n)$ denote the number of prime numbers no more than $n$. We have 
$$
\pi(n)=\frac{n}{\log n}+O(\frac{n}{\log^2n}).
$$ 
In particular, $\pi(n)\sim \frac{n}{\log n}$.
\end{lemma}
\bigskip

\begin{lemma}[Bertrand's postulate]\label{strong Bertrand}
For all $k\geq 89693$, there exists a prime number in the interval $$k<p\leq (1+\frac{1}{\log^3k})k.$$
\end{lemma}

\begin{proof}
The proof of this inequality appears in \cite{dusart2018explicit}.
\end{proof}

\begin{lemma}\label{prime number theorem 2}
Let $p_n$ denotes the $n^{th}$ prime number. We have 
$$
p_n=n\log n+O(n\log \log n).
$$ 
In particular, $p_n\sim n\log n.$
\end{lemma}
\bigskip

\begin{lemma}[The asymptotic squeeze principle]\label{special requirement}
There exist an $n_0\in \mathbb{N}$ such that for all even numbers $n\geq n_0$ there exist an axis
$$
\mathbb{L}_{[z_1],[z_2]}\inn \C^o(n,\CC_{\N})
$$ 
and 
$$
\mathbb{L}_{[w_1],[w_2]}\inn \C^o(n+t,\CC_{\N})
$$ 
with $\Re(z_1)\lesssim \Re(z_2)$ and $\Re(w_1)\lesssim \Re(w_2)$ such that $\Re(z_1)\lesssim \Re(w_1)$ and $\Re(z_2)\lesssim \Re(w_2)$ with $z_2\in \CC_{\p}$ and $w_1 \in \CC_{\p}$ for $t\geq 4$.
\end{lemma}

\begin{proof}
Let us set $\Re(z_2)$ to be a prime number and choose $\Re(z_2)$ to be the $\pi(\frac{3n}{4})^{th}$ prime number. We note that via the prime number theorem 
\begin{align}
\pi(\frac{3n}{4})&=\frac{\frac{3n}{4}}{\log (\frac{3n}{4})}+O(\frac{n}{\log^2n})\nonumber \\&=\frac{3n}{4\log n}+O(\frac{n}{\log^2n}).\nonumber
\end{align}
We also note that, via basic power series identities, we can write 
\begin{align}
-\log (1-\frac{1}{\log n})&=\frac{1}{\log n}+\frac{1}{2(\log n)^2}+\frac{1}{3(\log n)^3}+\cdots \nonumber \\&=\frac{1}{\log n}+O(\frac{1}{(\log n)^2}).\nonumber
\end{align}
By the lemma \ref{prime number theorem 2} and the lemma \ref{prime number inequality}, we obtain 
\begin{align}
\Re(z_2)=p_{\pi(\frac{3n}{4})}&=\pi(\frac{3n}{4})\log \pi(\frac{3n}{4})+O(\pi(\frac{3n}{4})\log \log \pi(\frac{3n}{4}))\nonumber \\&=(\frac{3n}{4\log n}+O(\frac{n}{\log^2n}))(\log (\frac{3n}{4\log n}+O(\frac{n}{\log^2n}))+O(\pi(\frac{3n}{4})\log \log \pi(\frac{3n}{4}))\nonumber 
\end{align}
We note that we can write 
\begin{align}
\log (\frac{3n}{4\log n}+O(\frac{n}{\log^2n})&=\log(\frac{3}{4})+\log(\frac{n}{\log n})+\log (1+O(\frac{1}{\log n}))\nonumber \\&=\log n-\log \log n+O(1)\label{main term 1}
\end{align}
From \eqref{main term 1}, we can write for the product 
\begin{align}
(\frac{3n}{4\log n}+O(\frac{n}{\log^2n}))(\log (\frac{3n}{4\log n}+O(\frac{n}{\log^2n}))&=(\frac{3n}{4\log n}+O(\frac{n}{\log^2n}))(\log n-\log \log n+O(1))\nonumber \\&=\frac{3n}{4}+O(\frac{n\log \log n}{\log n})
\end{align}
as the main term. Now, we analyze the error term in a similar manner. Using the prime number theorem, we can write 
\begin{align}
\pi(\frac{3n}{4})\log \log \pi(\frac{3n}{4})=(\frac{3n}{4\log n}+O(\frac{n}{\log ^2n}))(\log \log (\frac{3n}{4\log n}+O(\frac{n}{\log^2n}))).
\end{align}
We observe that 
\begin{align}
\log \log (\frac{3n}{4\log n}+O(\frac{n}{\log^2n}))&=\log( \log (\frac{3n}{4\log n})+\log (1+O(\frac{1}{\log n})))\ll \log \log n
\end{align}
so that we obtain for the product 
\begin{align}
\pi(\frac{3n}{4})\log \log \pi(\frac{3n}{4})&=(\frac{3n}{4\log n}+O(\frac{n}{\log ^2n}))(\log \log (\frac{3n}{4\log n}+O(\frac{n}{\log^2n})))\nonumber \\&\ll \frac{n\log \log n}{\log n}.\label{error term}
\end{align}
By combining \eqref{main term 1} and \eqref{error term}, we obtain $$
\Re(z_2)=\frac{3n}{4}+O(\frac{n\log \log n}{\log n})+O(\frac{n\log \log n}{\log n})=\frac{3n}{4}+O(\frac{n\log \log n}{\log n}).
$$
Consequently, we have for the real weight of the lower axis point 
\begin{align}
\Re(z_1)&=n-\Re(z_2)\nonumber\\&=n-\frac{3n}{4}+O(\frac{n\log \log n}{\log n})\nonumber \\&=\frac{n}{4}+O(\frac{n\log \log n}{\log n}).\label{inequality 2}
\end{align}
It is easy to see that 
$$
\Re(z_1)\sim \frac{n}{4}<\frac{n}{2}
$$
and 
$$
\Re(z_2)\sim \frac{3n}{4}>\frac{n}{2}
$$ 
Now, by virtue of the lemma \ref{strong Bertrand}, we set $\Re(w_1)$ to be a prime number and choose $\Re(w_1)$ so that
\begin{align}
\frac{n}{4}<\Re(w_1)&\leq (1+\frac{1}{\log^3 \frac{n}{4}})(\frac{n}{4})\label{inequality 3}
\end{align}
for all $n\geq 358772$. It implies that $\Re(z_1)\lesssim \Re(w_1)$. We deduce
$$
\Re(w_1)\lesssim \frac{n+t}{2}
$$
for $t\geq 4$, since 
\begin{align*}
(1+\frac{1}{\log^3\frac{n}{4}})(\frac{n}{4})\lesssim\frac{n}{2}
\end{align*} 
by virtue of the fact that 
$$
(1+\frac{1}{\log^3 \frac{n}{4}})\sim 1\quad (n\longrightarrow \infty).
$$
It follows from \eqref{inequality 3} the lower bound
\begin{align*}
\Re(w_2)&=n+t-\Re(w_1)\\
&\geq n+t-(1+\frac{1}{\log^3\frac{n}{4}})(\frac{n}{4})\\
&\mbox{and since }(1+\frac{1}{\log^3 \frac{n}{4}})\sim 1,~(n\longrightarrow \infty)\\
&\sim n-\frac{n}{4}+t\quad (n\longrightarrow \infty)\\
&> n-\frac{n}{4}=\frac{3n}{4}
\sim \Re(z_2)
\end{align*}
for all $t\geq 4$ and  $n>n_0$ for some fixed $n_0\in \mathbb{N}$. This completes the proof.
\end{proof}
\bigskip

We are now ready to prove the binary Goldbach conjecture for all even numbers greater than some $n_0\in \mathbb{N}$. This result provides an alternative solution to our first result and in very few instances adopts the proof technique in \cite{agama2022asymptotic}. The benefit of the strong version of the Bertrand postulate (Lemma \ref{strong Bertrand}) is good enough to verify the asymptotic version of the binary Goldbach conjecture using this version of the squeeze principle, which is a slight variation of the version that appears in the paper \cite{agama2022asymptotic}. 

\begin{theorem}[The asymptotic binary Goldbach theorem]\label{binary Goldbach theorem}
There exist some $n_0\in \mathbb{N}$ such that every even number $n\geq n_0$ can be written as a sum of two prime numbers.
\end{theorem}

\begin{proof}
We note that the above statement is equivalent to the statement that: the cCoPs $\C^o(n,\CC_{\p})$ are non--empty for all even $n\geq n_0$.

By the remark \ref{important remark} all cCoPs basing on $\CC_{\N}$ with generators $\geq 2$ are non--empty. Using the lemma \ref{special requirement}, all cCoPs $\C^o(n)$ and $\C^o(n+4)$ with even generators $n\geq n_0$ fulfill the requirements of the special squeeze principle (Lemma \ref{L-special_squeeze}) or the squeeze principle with $\mathbb{M}:=\mathbb{N}$ and $\mathbb{B}:=\mathbb{P}$. Hence, for each such $n$ there is always a non--empty cCoP $\C^o(n+2,\CC_{\p})$. We start with $\C^o(n_0)$ and $\C^o(n_0+4)$ and continue this procedure with $\C^o(n_0+k)$ and $\C^o(n_0+k+4)$ for all even $k\geq 2$. We verify that all cCoPs $\C^o(n_0+k+2)$ for even $k\geq 2$ ad infinitum are non--empty. 
\end{proof}
\bigskip

\section{The generalized squeeze principle and applications}\label{sec:generalized}

In this section, we generalize the \emph{squeeze principle} introduced and developed in the previous section.

\begin{definition}
Let $\mathbb{A}\subset \mathbb{N}$. We call the set 
$$
\mathcal{C}(n,\bigotimes_{i=1}^{h}\mathbb{A}_i):=\left \{[x_1],[x_2],\ldots,[x_h]~\bigg |~x_i\in \mathbb{A}_i,~n=\sum \limits_{i=1}^{h}x_i\right\}
$$ 
a \emph{multivariate circle of partition} generated by $n\in \mathbb{N}$ with base \emph{regulators} $\bigotimes_{i=1}^{h}\mathbb{A}_i$ the $h$-fold direct product of the sets $\mathbb{A}_i$. We call members of the multivariate circle of partitions \emph{multivariate points}. We denote the \emph{weight} of each point by $||[x_i]||:=x_i\in \mathbb{A}_i$ and for the corresponding weight set of the multivariate circle of partitions 
$$
||\mathcal{C}(n,\bigotimes_{i=1}^{h}\mathbb{A}_i)||:=\{(x_1,x_2,\ldots x_h)\in \bigotimes_{i=1}^{h}\mathbb{A}_i~|~\sum \limits_{i=1}^{h}x_i=n\}.
$$
\bigskip

We denote an \emph{axis} of the multivariate circle of partition $\mathcal{C}(n,\bigotimes_{i=1}^{h}\mathbb{A}_i)$ by $\mathbb{L}_{[x_1],[x_2],\ldots,[x_h]}$ if and only if  $x_i\in \mathbb{A}_i$ for each $1\leq i\leq h$ and 
$$
n=\sum \limits_{i=1}^{h}x_i.
$$
We say that the axis points $[x_i]$ for each $1\leq i\leq h$ are axis \emph{residents}. We do not view the axes to be different among other axis of the form $\mathbb{L}_{[x_1],[x_2],\ldots, [x_h]}$ up to the rearrangements of its residents points. In special cases where the points 
$$
[x_k] \in \mathcal{C}(n,\bigotimes_{i=1}^{h}\mathbb{A}_i)
$$ 
is such that $hx_k=n$, we call $[x_i]$ the \emph{center} of the multivariate circle of partitions. If it exists, then we call it a \emph{degenerated axis} $\mathbb{L}_{[x_k]}$ in comparison with the \emph{real axes} $\mathbb{L}_{[x_1],[x_2],\ldots, [x_h]}$, where not all of the weights $x_i$ can be equal. We denote the assignment of an axis $\mathbb{L}_{[x_1],[x_2],\ldots, [x_h]}$ to the multivariate CoP $\mathcal{C}(n,\bigotimes_{i=1}^{h}\mathbb{A}_i)$ by
$$
\mathbb{L}_{[x_1],[x_2],\ldots, [x_h]} \inn \mathcal{C}(n,\bigotimes_{i=1}^{h}\mathbb{A}_i)
$$
which means
$$
[x_1],[x_2],\ldots,[x_h] \in \mathcal{C}(n,\bigotimes_{i=1}^{h}\mathbb{A}_i)
$$
with 
$$
n=\sum \limits_{i=1}^{h}x_i
$$
for a fixed $n\in \mathbb{N}$ with $x_i\in \mathbb{A}_i$ for each $1\leq i\leq h$ or vice versa. We denote the number of real axes of the multivariate circle of partitions by
$$
\nu(n,\bigotimes_{i=1}^{h}\mathbb{A}_i):=\#\lbrace \mathbb{L}_{[x_1],[x_2],\ldots, [x_h]} \inn  \mathcal{C}(n,\bigotimes_{i=1}^{h}\mathbb{A}_i) \mid x_i\neq x_j\rbrace
$$
for all $1<i<j\leq h$.
\bigskip

 In the special case where we fix $h=2$ and take $\mathbb{A}_i=\mathbb{A}\subset \mathbb{N}$, we obtain the circle of partitions 
 $$
 \mathcal{C}(n,\mathbb{A}):=\left \{[x]~|~x,n-x\in \mathbb{A}\right \}
 $$ 
 and the corresponding counting function for the axes set 
 $$
 \nu(n,\mathbb{A}):=\#\lbrace \mathbb{L}_{[x],[y]}\inn\mathcal{C}(n,\mathbb{A})\mid x\neq y\rbrace.
 $$
 \end{definition}
 \bigskip
 
 This structure was studied extensively in \cite{CoP} in the case where one allows only two axis points on their axes. The squeeze principle is the statement 

\begin{theorem}[The squeeze principle]
Let $\mathbb{B} \subset\mathbb{M}\subseteq\mathbb{N}$ and
$\mathcal{C}(m,\mathbb{B})$ and $\mathcal{C}(m+t,\mathbb{B})\neq\emptyset$
for $t\geq 4$. 
If there exists $\mathbb{L}_{[x],[y]}\inn \mathcal{C}(m+t,\mathbb{M})$ with $x\in \mathbb{B}$ and $x<y$ such that 
\begin{align}\label{E_maxcovered}
y>w:= \mathrm{\max}\{u\in||\mathcal{C}(m,\mathbb{M})||\mid u\in \mathbb{B}\}>m-x,
\end{align}
then there exists $\mathcal{C}(s,\mathbb{B})\neq\emptyset$
such that $m<s<m+t$.
\end{theorem}

Indeed, it has found some unexpected applications to \textbf{Goldbach}-type problems and has been applied to study additive prime number problems requiring certain partitions into certain subsets of the positive integers. The power of this principle allows one to exhaust any interval of the form $[n,n+t]$ for $t\geq 4$ in way to conclude that all even numbers in this interval can be written as the sum of two prime numbers. In general, the set may be extended to a general subset of the integers and it mostly suffices to check if the underlying conditions of the principle are all satisfied in order to run this test. In the previous section, we proved:

\begin{lemma}[The squeeze principle]
Let $\mathbb{B}\subset\mathbb{M}\subseteq \mathbb{N}$ and
$\mathcal{C}^o(n,\mathbb{C}_{\mathbb{M}})$ and $\mathcal{C}^o(n+t,\mathbb{C}_{\mathbb{M}})$ with $t\geq 4$ be non--empty cCoPs with integers $n,t,s$ of the same parity.
If there exist an axis $\mathbb{L}_{[v_1],[w_1]}\inn\mathcal{C}^o(n,\mathbb{C}_{\mathbb{M}})$ with $w_1\in \mathbb{C}_{\mathbb{B}}$ and an axis $\mathbb{L}_{[v_2],[w_2]}\inn \mathcal{C}^o(n+t,\mathbb{C}_{\mathbb{M}})$ with $v_2\in \mathbb{C}_{\mathbb{B}}$ such that 
\begin{align}
\Re(v_1)<\Re(v_2)\quad \text{and}\quad \Re(w_1)<\Re(w_2)\nonumber
\end{align}
then there exists an axis $\mathbb{L}_{[\Re(v_2)],\Re([w_1)]}\inn \mathcal{C}(n+s,\mathbb{B})$ with $0<s<t$. Hence $\mathcal{C}^o(n+s,\mathbb{C}_{\mathbb{M}})$ is also not empty.
\end{lemma}

We obtain an analogous versions of the \emph{squeeze} principle in the setting of axes of circle of partitions with at least two resident points.

\begin{lemma}[Generalized squeeze principle]\label{general squeeze}
Let $\mathbb{A}\subset \mathbb{N}$ with $\mathcal{C}(n,\bigotimes_{i=1}^{h}\mathbb{A})\neq \emptyset$ for a fixed $n\in \mathbb{H}\subseteq \mathbb{N}.$ If $t\in \mathbb{N}$ is such that $n$ and $n+t$ are not consecutive integers in $\mathbb{H}$ and there exists an axes 
$$
\mathbb{L}_{[x_1],[x_2],\ldots,[x_h]}\inn \mathcal{C}(n+t,\bigotimes_{i=1}^{h}\mathbb{N})
$$
with $x_i\in \mathbb{A}$ for all $1\leq i\leq h-1$ and $x_i<x_h$ for all $1\leq i\leq h-1$ such that 
$$
x_h>w:=\mathrm{\max}\{u\in ||\mathcal{C}(n,\mathbb{N})||~|~u\in \mathbb{A}\}>n-\sum \limits_{i=1}^{h-1}x_i
$$ 
then there exists 
$$
\mathcal{C}(s,\bigotimes_{i=1}^{h}\mathbb{A})\neq\emptyset
$$ 
for $n<s<n+t$ with $s\in \mathbb{H}$.
\end{lemma}

\begin{proof}
We note that from the hypothesis 
$$
x_h>w:=\mathrm{\max}\{u\in ||\mathcal{C}(n,\mathbb{N})||~|~u\in \mathbb{A}\}>n-\sum \limits_{i=1}^{h-1}x_i
$$ 
we can write 
$$
n=w+(n-w)<w+\sum \limits_{i=1}^{h-1}x_i<\sum \limits_{i=1}^{h}x_i=n+t
$$ 
for  $x_i<x_h$ for all $1\leq i\leq h-1$. Clearly $w\in \mathbb{A}$ with $x_i \in \mathbb{A}$ for $1\leq i\leq h-1$ so that there exists an axis 
$$
\mathbb{L}_{[w],[x_1],\ldots,[x_{h-1}]}\inn \mathcal{C}(s,\bigotimes_{i=1}^{h}\mathbb{A})
$$ 
with $n<s<n+t$ and $s\in \mathbb{H}$, since $n,n+t$ are not consecutive in $\mathbb{H}$. This means that $s$ can be written as the sum of $h$~(not all possibly distinct) elements of $\mathbb{A}\subset \mathbb{N}$.
\end{proof}
\bigskip

We show how this principle can be applied to solve \textbf{Goldbach}-type problems with $h$ summands for $h\geq 2$.

\begin{theorem}[The partition law]
Let $\mathbb{A}\subset \mathbb{N}$ and suppose that $\mathcal{C}(n,\bigotimes_{i=1}^{h}\mathbb{A})\neq \emptyset$ for infinitely many $n\in \mathbb{H}\subseteq \mathbb{N}$. If for each $t\in \mathbb{N}$ such that $n,n+t$ are not consecutive in $\mathbb{H}$ there exists at least an axis 
$$
\mathbb{L}_{[x_1],[x_2],\ldots,[x_h]}\inn \mathcal{C}(n+t,\bigotimes_{i=1}^{h}\mathbb{N})
$$ 
with $x_i\in \mathbb{A}$ for all $1\leq i\leq h-1$ and $x_i<x_h$ for all $1\leq i\leq h-1$ such that 
$$
x_h>w:=\mathrm{\max}\{u\in ||\mathcal{C}(n,\mathbb{N})||~|~u\in \mathbb{A}\}>n-\sum \limits_{i=1}^{h-1}x_i
$$
then there are multivariate circles of partitions with the property $$
\mathcal{C}(s,\bigotimes_{i=1}^{h}\mathbb{A})\neq \emptyset
$$ 
for all $s\geq k$ for a fixed $k\in \mathbb{N}$ with $s\in \mathbb{H}$, which means every number in $\mathbb{H}\subset \mathbb{N}$ $\geq k$ can be written as the sum of $h$ elements~(not all possibly distinct)~of $\mathbb{A}$.
\end{theorem}

\begin{proof}
Suppose that $\mathbb{A}\subset \mathbb{N}$ and let $k$ be the smallest number in $\mathbb{H}\subseteq \mathbb{N}$ such that $\mathcal{C}(k,\bigotimes_{i=1}^{h}\mathbb{A})\neq \emptyset$, which means $k$ is the smallest number in $\mathbb{H}$ such that it can be written as the sum of $h$ elements in $\mathbb{A}\subset \mathbb{N}$. Let us choose $t_0\in \mathbb{N}$ such that $k,k+t_0$ are not consecutive in $\mathbb{H}$. By the hypothesis there exists at least an axis 
$$
\mathbb{L}_{[x_1],[x_2],\ldots,[x_h]}\inn \mathcal{C}(k+t_0,\bigotimes_{i=1}^{h}\mathbb{N})
$$ 
with $x_i\in \mathbb{A}$ for all $1\leq i\leq h-1$ and $x_i<x_h$ for all $1\leq i\leq h-1$ such that
$$
x_h>w:=\mathrm{\max}\{u\in ||\mathcal{C}(k,\mathbb{N})||~|~u\in \mathbb{A}\}>k-\sum \limits_{i=1}^{h-1}x_i.
$$ 
It follows from the lemma \ref{general squeeze} that there exists an $s\in \mathbb{N}$ with $k<s<k+t_0$ such that 
$$
\mathcal{C}(s,\bigotimes_{i=1}^{h}\mathbb{A})\neq \emptyset
$$
which means $s\in \mathbb{H}$ can be written as the sum of $h$~(not all) possibly distinct elements of $\mathbb{A}$. We consider the sub-intervals $[k,s]$ and $[s,k+t_0]$. If $k<s:=k+t_1$ are not consecutive in $\mathbb{H}$, then there exists at least an axis 
$$
\mathbb{L}_{[x_1],[x_2],\ldots,[x_h]}\inn \mathcal{C}(k+t_1,\bigotimes_{i=1}^{h}\mathbb{N})
$$ 
with $x_i\in \mathbb{A}$ for all $1\leq i\leq h-1$ and $x_i<x_h$ for all $1\leq i\leq h-1$ such that 
$$
x_h>w:=\mathrm{\max}\{u\in ||\mathcal{C}(k,\mathbb{N})||~|~u\in \mathbb{A}\}>k-\sum \limits_{i=1}^{h-1}x_i.
$$ 
It follows from the lemma \ref{general squeeze} that there exists some $u\in \mathbb{H}$ with $k<u<k+t_1$ such that 
$$
\mathcal{C}(u,\bigotimes_{i=1}^{h}\mathbb{A})\neq \emptyset
$$ 
which means $u\in \mathbb{H}$ can also be written as the sum of $h$~(not all) possibly distinct elements of $\mathbb{A}$. We can similarly iterate this argument on the intervals $[k,u],[u,s],[s,k+t_0]$ in so far as there exists an element of $\mathbb{H}$ in any of the intervals. By virtue of the hypothesis, the $t\in \mathbb{N}$ can be arbitrarily chosen such that $k,k+t$ are not consecutive in $\mathbb{H}$ and with the existence of at least an axis 
$$
\mathbb{L}_{[x_1],[x_2],\ldots,[x_h]}\inn \mathcal{C}(k+t_0,\bigotimes_{i=1}^{h}\mathbb{N})
$$ 
with $x_i\in \mathbb{A}$ for all $1\leq i\leq h-1$ and $x_i<x_h$ for all $1\leq i\leq h-1$ such that 
$$
x_h>w:=\mathrm{\max}\{u\in ||\mathcal{C}(k,\mathbb{N})||~|~u\in \mathbb{A}\}>k-\sum \limits_{i=1}^{h-1}x_i.
$$ 
The iterative arguments can be extended to all elements of $\mathbb{H}$ under the assumption that 
$$
\mathcal{C}(n,\bigotimes_{i=1}^{h}\mathbb{A})\neq \emptyset
$$ 
for infinitely many $n\in \mathbb{H}\subseteq \mathbb{N}$. This means that every element $n\in \mathbb{H}$ can be written as the sum of $h$ elements of $\mathbb{A}$, not all distinct.
\end{proof}
%%%%%%%%%%%%%%%%%%%%%%%%%%%%%%%%%%%%%%%%%%%%%%%%%%%%%%%%%%%%%%%%%%%%%%%%
\rule{100pt}{1pt}

\bibliographystyle{amsplain}

\end{document}